\documentclass[a4paper,11pt]{article}
\usepackage{graphics}
\usepackage{amssymb}
\usepackage{amsmath}
\usepackage{anysize}
\usepackage{graphics}

\begin{document}
\begin{center}
\textsc{\Large The Generating Function of the Embedding Capacity for 4-dimensional Symplectic Ellipsoids}\vspace{5mm}
\\{\large David Bauer\footnote[1]{Max Planck Institute for Mathematics in the Sciences, Leipzig,
Germany; \textsf{dbauer@mis.mpg.de}}}
\end{center}\vspace{2mm}

\newcounter{sect}\setcounter{sect}{0}
\newcounter{def}\setcounter{def}{0}
\newcounter{lemma}\setcounter{lemma}{0}
\newcounter{theorem}\setcounter{theorem}{0}
\newcounter{problem}\setcounter{problem}{0}
\newcounter{corollary}\setcounter{corollary}{0}
\newcounter{proposition}\setcounter{proposition}{0}
\newcommand{\writesection }[1]{\stepcounter{sect}\textbf{\normalsize\thesect .\hspace{1mm} #1.\ }}
\newcommand{\writedefinition }{\stepcounter{def}\textsc{\normalsize\\ \hspace{6mm}Definition \thedef .\hspace{2mm}}}
\newcommand{\writelemma }{\stepcounter{lemma}\textsc{\normalsize\\ \hspace{6mm}Lemma \thelemma .\hspace{2mm}}}
\newcommand{\writetheorem }{\stepcounter{theorem}\textsc{\normalsize\\ \hspace{6mm}Theorem \thetheorem .\hspace{2mm}}}
\newcommand{\writeproblem }{\stepcounter{problem}\textsc{\normalsize\\ \hspace{6mm}Problem \theproblem .\hspace{2mm}}}
\newcommand{\writecorollary }{\stepcounter{corollary}\textsc{\normalsize\\ \hspace{6mm}Corollary \thecorollary
.\hspace{2mm}}}
\newcommand{\writeproposition }{\stepcounter{proposition}\textsc{\normalsize\\ \hspace{6mm}Proposition \theproposition
.\hspace{2mm}}}

\numberwithin{equation}{sect}
\renewcommand{\labelenumi}{(\alph{enumi})}

\small

\begin{abstract}
\footnotesize Quite recently, McDuff showed that the existence of a symplectic embedding of one four-dimensional
ellipsoid into another can be established by comparing their corresponding sequences of ECH capacities. In this note we
show that these sequences can be encoded in a generating function, which gives several new equivalent formulations of
McDuff's theorem.
\end{abstract}\vspace{3mm}

\writesection{Embedding 4-dimensional Symplectic Ellipsoids} We consider ellipsoids
\begin{align*}
E(a,b):=\left\{z\in\mathbb{C}^2:\frac{|z_1|^2}{a}+\frac{|z_2|^2}{b}\leq 1\right\}
\end{align*}
equipped with the standard symplectic structure $\omega_0=\mathrm{d}x_1\wedge\mathrm{d}y_1+\mathrm{d}x_2\wedge\mathrm{d}y_2$ of Euclidean space $\mathbb{R}^4$. The embedding problem in symplectic geometry asks if for given integers $a,b,c,d>0$ there exists a symplectic embedding
$\mathrm{int}\,E(a,b)\stackrel{s}{\hookrightarrow} E(c,d)$. Since each such embedding preserves the volume, an immediate obstruction for existence is $ab\leq cd$. \vspace{1mm}

There are further obstructions which have their origin in embedded contact homology. Namely, define $\mathcal{N}(a,b)$ to be the sequence of numbers from the set
\begin{align*}
\mathcal{S}(a,b):=\{ka+lb:k,l\in\mathbb{Z}\mbox{ and }k,l\geq 0\}
\end{align*}
arranged in nondecreasing order with repetitions. For example, we have
\begin{align*}
\mathcal{N}(2,3)=(0,2,3,4,5,6,6,7,8,8,9,9,\ldots).
\end{align*}
For sequences of numbers $\mathcal{A}$ and $\mathcal{B}$ define a partial ordering by saying
$\mathcal{A}\preceq\mathcal{B}$ if, for all $n\geq 0$, the $n$-th entry of $\mathcal{A}$ is not larger than the $n$-th
entry of $\mathcal{B}$. Hutchings showed in \cite{Hutc2} that an obstruction for the embedding problem is given by
$\mathcal{N}(a,b)\preceq\mathcal{N}(c,d)$. Indeed, as conjectured by Hofer and recently proved by McDuff in
\cite{Duff2}, this is the only obstruction.

\writetheorem\textsl{There is a symplectic embedding $\mathrm{int}\,E(a,b)\stackrel{s}{\hookrightarrow} E(c,d)$ if and only if
\begin{align*}
\mathcal{N}(a,b)\preceq\mathcal{N}(c,d).
\end{align*}
}\vspace{2mm}

\begin{figure}[htbp]
\begin{minipage}{40mm}
\hspace{40mm}
\end{minipage}
\begin{minipage}{60mm}
\includegraphics[5cm,4cm]{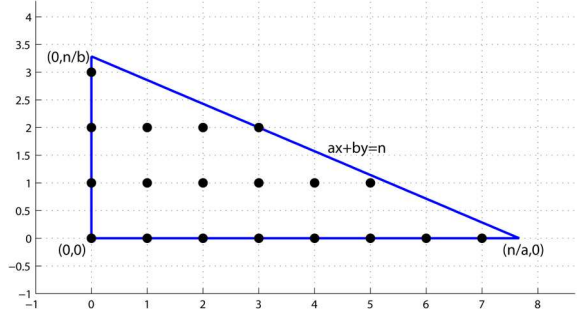}
\end{minipage}
\caption{Interpreting $\mathcal{L}_n(a,b)$ as a lattice count.}\label{fig:TriLatPts}
\end{figure}

Hence the embedding problem for symplectic ellipsoids can be reduced to studying the sequences $\mathcal{N}(a,b)$.
Define a new sequence $\mathcal{L}(a,b)$ by
\begin{align*}
\mathcal{L}_n(a,b):=\max\{j:\mathcal{N}_j(a,b)\leq n\}=\#\{m\in\mathcal{S}(a,b):m\leq n\}.
\end{align*}
From the definition it is clear, that
\begin{align}\label{eq.EquivSeqLN}
\mathcal{L}(a,b)\succeq\mathcal{L}(c,d) \Longleftrightarrow \mathcal{N}(a,b)\preceq\mathcal{N}(c,d).
\end{align}
Geometrically, $\mathcal{L}_n(a,b)$ corresponds to the number of lattice points in the triangle $T_{a,b}^n$ bounded by
$x=0$, $y=0$ and $ax+by=n$, including points on its boundary (Figure \ref{fig:TriLatPts}). \vspace{1mm}

The aim of this note is to remark that the generating function of $\mathcal{L}(a,b)$ is given by a surprisingly simple
formula.

\writeproposition\textsl{For $0\leq z<1$ we have the expansion
\begin{align}\label{eq.PowSerL}
\frac{1}{(1-z)(1-z^a)(1-z^b)}=\sum\limits_{n=0}^\infty \mathcal{L}_n(a,b)z^n.
\end{align}
}\vspace{1mm}

{\sc Proof:} We have
\begin{align*}
\frac{1}{(1-z)(1-z^a)(1-z^b)}&=\left(\sum\limits_{k=0}^\infty z^k\right)\left(\sum\limits_{l=0}^\infty
z^{al}\right)\left(\sum\limits_{m=0}^\infty z^{bm}\right)\\&= \sum\limits_{n=0}^\infty
\left(\#\left\{(k,l,m)\in\mathbb{Z}^3: k,l,m\geq 0\mbox{ and } k+al+bm=n\right\}\right)z^n\\&= \sum\limits_{n=0}^\infty
\left(\#\left\{(l,m)\in\mathbb{Z}^2: l,m\geq 0\mbox{ and } al+bm\leq n\right\}\right)z^n=\sum\limits_{n=0}^\infty
\mathcal{L}_n(a,b)z^n.
\end{align*} \hfill
$\Box$\vspace{4mm}

There is also a geometric interpretation behind this formula, which will be explained in the next section. Note that
$\mathcal{L}_n(a,b)$ corresponds to the number of partitions of $n$ into parts of size $1$, $a$ or $b$ which is known
as a denumerant problem. In this case one always obtains a rational generating function with poles that are roots of
unity. Multiplying both sides of (\ref{eq.PowSerL}) by the denominator and comparing coefficients leads to the linear
recurrence relation
\begin{align*}
\mathcal{L}_{n}(a,b)=&\mathcal{L}_{n-1}(a,b)+\mathcal{L}_{n-a}(a,b)+\mathcal{L}_{n-b}(a,b)+\mathcal{L}_{n-a-b-1}(a,
b)\\&-\mathcal{L}_{n-a-1}(a,b)-\mathcal{L}_{n-b-1}(a,b)-\mathcal{L}_{n-a-b}(a,b)
\end{align*}
for $n>0$. To initiate we take $\mathcal{L}_0(a,b)=1$ and set $\mathcal{L}_n(a,b):=0$ for $n<0$.
The following relation can be proved in an elementary way (see \cite{Hard}, section 5.6).

\writeproposition\textsl{For $n>0$ we have
\begin{align}\label{eq.RecCoeffL}
\mathcal{L}_n(a,b)=\mathcal{L}_{n-1}(a,b)+\left\lfloor\frac{n}{ab}\right\rfloor+\varepsilon(n)
\end{align}
where $\varepsilon(n)$ is either $0$ or $1$ and its value just depends on the remainder
\begin{align*}
[n]\in\frac{\mathbb{Z}}{ab\mathbb{Z}}.
\end{align*}}\vspace{2mm}

In some sense the whole information of $\mathcal{L}(a,b)$ is therefore stored in its first $ab$ terms. Moreover, one
obtains the asymptotic behaviour
\begin{align*}
\mathcal{L}_n(a,b)\sim\frac{n^2}{2ab}.
\end{align*}
In the following, we denote the generating function by
\begin{align*}
 g_{a,b}(z)=\frac{1}{(1-z)(1-z^a)(1-z^b)}.
\end{align*}
Denote further by $f^{(k)}$ the $k$-th derivative of a function $f$. Via Cauchy's integral formula we compute
\begin{align*}
\mathcal{L}_n(a,b)=\frac{g_{a,b}^{(n)}(0)}{n!}=\frac{1}{2\pi\mathrm{i}}\int_\gamma\frac{g_{a,b}(\xi)\,\mathrm{d}\xi}{
\xi^ { n+1 }
}=\frac{1}{2\pi\mathrm{i}}\int_\gamma\frac{\mathrm{d}\xi}{(1-\xi)(1-\xi^a)(1-\xi^b)\xi^{n+1}},
\end{align*}
which might be useful for numerical purposes.\vspace{1mm}

On the space $\mathcal{C}^\infty((-1,1),\mathbb{R})$ consider the partial ordering by saying $f\preceq g$ iff
$f^{(k)}(x)\leq g^{(k)}(x)$ for all $k\geq 0$ and $x\in[0,1)$. Putting things together we obtain the following

\writecorollary\textsl{There is a symplectic embedding $\mathrm{int}\,E(a,b)\stackrel{s}{\hookrightarrow} E(c,d)$ if and
only if one of the following equivalent conditions is fulfilled:}
\begin{enumerate}
\item $\mathcal{N}(a,b)\preceq\mathcal{N}(c,d)$
\item $\mathcal{L}(a,b)\succeq\mathcal{L}(c,d)$
\item $g_{a,b}\succeq g_{c,d}$
\end{enumerate}
\vspace{1mm}

{\sc Proof:} The equivalence of (a) and (b) was already noticed in ($\ref{eq.EquivSeqLN}$). Now (b) implies for
any integer $k\geq 0$ and $z\in[0,1)$
\begin{align*}
g_{a,b}^{(k)}(z)=\sum\limits_{n=k}^\infty k!\binom n k\mathcal{L}_n(a,b)z^{n-k}\geq\sum\limits_{n=k}^\infty k!\binom n
k\mathcal{L}_n(c,d)z^{n-k}=g_{c,d}^{(k)}(z).
\end{align*}
On the other hand (c) leads to
\begin{align*}
\mathcal{L}_k(a,b)=\frac{g_{a,b}^{(k)}(0)}{k!}\geq\frac{g_{c,d}^{(k)}(0)}{k!}=\mathcal{L}_k(c,d).
\end{align*}\hfill
$\Box$\vspace{4mm}

Thus the embedding question $\mathrm{int}\,E(a,b)\stackrel{s}{\hookrightarrow} E(c,d)$ relates to the problem if all
coefficients of
\begin{align*}
G_{a,b,c,d}(z):=\frac{(1-z^c)(1-z^d)-(1-z^a)(1-z^b)}{(1-z)(1-z^a)(1-z^b)(1-z^c)(1-z^d)}=g_{a,b}(z)-g_{c,d}
(z)=\sum\limits_ { n=0 } ^\infty \left(\mathcal{L}_n(a,b)-\mathcal{L}_n(c,d)\right)z^n
\end{align*}
are nonnegative. Since $G_{a,b,c,d}$ is again a rational function, its coefficients satisfy a linear recurrence. In
\cite{Gerh}, Conjecture 2 it is conjectured that each rational function, whose dominating poles (i.e. the ones of
maximal modulus) do not lie on $\mathbb{R_+}$, has infinitely many positive and infinitely many negative coefficients
in its power series expansion. Of course, we cannot apply this to $G_{a,b,c,d}$, since all of its poles have modulus 1
and $1\in\mathbb{R_+}$ occurs among them. One of the most celebrated results in the theory of linear recurrence
sequences is the Skolem-Mahler-Lech theorem. It asserts that if a sequence $(a_n)$ satisfies a linear recurrence
relation, then the zero set
\begin{align*}
\{n\in\mathbb{N}:a_n=0\}
\end{align*}
is the union of a finite set and finitely many arithmetic progressions.\vspace{1mm}

Let us use the approach via generating functions to check algebraically that for each positive integer $n\in\mathbb{N}$
there is a symplectic embedding 
\begin{align*}
\mathrm{int}\,E(1,n^2)\stackrel{s}{\hookrightarrow} B(n).
\end{align*}
Here the latter denotes the ball $B(n):=E(n,n)$ of radius $n$. Geometrically, this corresponds to a filling of $B(n)$ by
$n^2$ equal symplectic balls (Proposition 2.2 in \cite{Duff1}). The possibility of such a filling can be quite easily
observed via toric models. For details we refer the reader to the survey paper \cite{Duff1}.\vspace{1mm}

With the lattice count interpretation we have
\begin{align*}
\mathcal{L}_k(n,n)=d\left(\left\lfloor\frac{k}{n}\right\rfloor\right),
\end{align*}
where $d(k):=\frac{1}{2}(k+1)(k+2)$ denotes the $k$-th triangle number. Consequently, by Proposition 1
\begin{align*}
g_{n,n}(z)=\frac{1}{(1-z)(1-z^n)^2}=\sum\limits_{k=0}^\infty d\left(\left\lfloor\frac{k}{n}\right\rfloor\right)z^k.
\end{align*}
For integers $k\geq 0$ set
\begin{align*}
c(k)=\Bigg\{\begin{array}{rl}
                 \hspace{3mm} 1 & \mbox{ if }k\equiv 0\ (\mathrm{mod}\ n),\\
                  -1 & \mbox{ if }k\equiv 1\ (\mathrm{mod}\ n), \\
                 0 & \mbox{ otherwise.}
               \end{array}
\end{align*}
Then
\begin{align*}
\frac{(1-z^n)^2}{(1-z)(1-z^{n^2})}&=\frac{1-z^n}{1-z}\cdot(1-z^n)\sum\limits_{k=0}^\infty z^{kn^2}= \left(1+z+\ldots +z^{n-1}\right)\sum\limits_{k=0}^\infty \left(z^{kn^2}-z^{(kn+1)n}\right) \\&=\sum\limits_{k=0}^\infty c\left(\left\lfloor\frac{k}{n}\right\rfloor\right)z^k,
\end{align*}
such that
\begin{align*}
g_{1,n^2}(z)=\frac{g_{1,n^2}(z)}{g_{n,n}(z)}\cdot g_{n,n}(z)=\frac{(1-z^n)^2}{(1-z)(1-z^{n^2})}\cdot g_{n,n}(z)=\left(\sum\limits_{k=0}^\infty c\left(\left\lfloor\frac{k}{n}\right\rfloor\right)z^k\right)\left(\sum\limits_{l=0}^\infty d\left(\left\lfloor\frac{l}{n}\right\rfloor\right)z^l\right).
\end{align*}
In view of (\ref{eq.PowSerL}) it suffices to show for each nonnegative integer $N$
\begin{align}\label{eq.convSeqcd}
\sum\limits_{k=0}^N c\left(\left\lfloor\frac{k}{n}\right\rfloor\right)d\left(\left\lfloor\frac{N-k}{n}\right\rfloor\right)\geq d\left(\left\lfloor\frac{N}{n}\right\rfloor\right).
\end{align}
For given $N\geq 0$ we pick integers $0\leq p,q,r$ with $q,r<n$ such that $N=pn^2+qn+r$. Setting $d(-1)=d(-2):=0$, we obtain from the periodicity of $c(k)$
\begin{align*}
\sum\limits_{k=0}^N
c\left(\left\lfloor\frac{k}{n}\right\rfloor\right)d\left(\left\lfloor\frac{N-k}{n}\right\rfloor\right)=&
\sum\limits_{j=0}^p\left((r+1)d(jn+q)+(n-r-1)d(jn+q-1)\right)
\\&-\sum_{j=0}^p\left((r+1)d(jn+q-1)+(n-r-1)d(jn+q-2)\right)\\=&\sum\limits_{j=0}
^p\left((r+1)(jn+q+1)+(n-r-1)(jn+q)\right) \\=&\frac{p(p+1)}{2}n^2+(p+1)qn+(p+1)(r+1)=(p+1)(N+1)-\frac{p(p+1)}{2}n^2.
\end{align*}
For $q<n$, $n\geq 2$ we have
\begin{align*}
\frac{3q}{2}+\frac{q^2}{2}=\frac{q(q+1)}{2}+q\leq\frac{nq}{2}+\frac{nq}{2},
\end{align*}
such that $\frac{3q}{2}+\frac{q^2}{2}\leq qn$ holds for all nonnegative integers $q<n$. One also easily checks that
$\frac{3pn}{2}\leq\frac{pn^2}{2}+p$ holds for all nonnegative integers $p,n$. Thus
\begin{align*}
(p+1)(N+1)\geq& (p+1)(pn^2+qn+1)=p^2n^2+pn^2+pqn+qn+p+1\\ \geq&
p^2n^2+\frac{pn^2}{2}+\frac{3pn}{2}+pqn+\frac{q^2}{2}+\frac{3q}{2}+1=\frac{(pn+q+1)(pn+q+2)}{2}+\frac{p(p+1)}{2}n^2\\=&
d\left(\left\lfloor\frac{N}{n}\right\rfloor\right)+\frac{p(p+1)}{2}n^2
\end{align*}
shows that (\ref{eq.convSeqcd}) is valid.\vspace{1mm}

The symplectic capacity function $c:[1,\infty)\rightarrow\mathbb{R}$ defined by
\begin{align*}
c(a):=\inf\left\{\mu:\mathrm{int}\,E(1,a)\stackrel{s}{\hookrightarrow} B(\mu)\right\}
\end{align*}
is studied in detail in \cite{DuffSchl}. We just computed $c(a^2)=a$ for positive integers $a$. Indeed, $c(a)=\sqrt{a}$ holds for $a\in\mathbb{N}$ if $a$ is 1,4 or $\geq 9$. The other values for integral $a$ are given by
\begin{align*}
c(2)=c(3)=2,\ c(5)=c(6)=\frac{5}{2},\ c(7)=\frac{8}{3},\ c(8)=\frac{17}{6}.
\end{align*}
We finish this section by remarking that Theorem 1 does not hold in higher dimensions. Counterexamples are due to Guth
\cite{Guth} and Hind-Kerman \cite{Hind}. Even worse, embedded contact homology only exists in dimension 4 and
there is so far no good guess of what a criterion for embedding ellipsoids could be.
\vspace{6mm}

\writesection{Counting Lattice Points in Polyhedra} Let $P\subset\mathbb{R}^d$ be a polyhedron. In order to count the lattice points in $P$ one associates the generating function
\begin{align*}
\sum\limits_{m\in P\cap\mathbb{Z}^d}\mathbf{x}^m\hspace{3mm}\mbox{with }\mathbf{x}^m=x_1^{\mu_1}\cdots x_d^{\mu_d}
\end{align*}
for $m=(\mu_1,\ldots ,\mu_d)$. The total number of lattice points in $P$ is then given by the value of the generating function at $\mathbf{x}=(1,\ldots ,1)$. The advantage of this approach is that these generating functions can still be computed for cones $K\subset\mathbb{R}^d$, which actually contain an infinite number of lattice points. A cone is characterized by the property that $0\in K$ and for every $x\in K$ and $\lambda\geq 0$ one has $\lambda x\in K$. For example, the generating function of the non-negative orthant is given by
\begin{align*}
\sum\limits_{m\in\mathbb{R}_+^d\cap\mathbb{Z}^d}\mathbf{x}^m=\prod\limits_{i=1}^d\frac{1}{1-x_i}.
\end{align*}
The generating function of a polyhedron $P$ is calculated as the sum of generating functions of tangent cones at the vertices of $P$, for details see \cite{Barv}.\vspace{1mm}

Usually a cone $K$ is given as a span of vectors $u_1,\ldots ,u_k\in\mathbb{R}^d$,
\begin{align*}
K=\mathrm{co}(u_1,\ldots ,u_k),
\end{align*}
meaning that every vector $v\in K$ can be written as a sum $v=\sum\lambda_iv_i$ with $\lambda_i\geq 0$. A cone $K$ is called unimodular, if it is spanned by $u_1,\ldots ,u_d\in\mathbb{Z}^d$ and these vectors form a basis of the lattice. Generating functions for unimodular cones are particularly easy to calculate. Unfortunately, all tangent cones of the triangle $T_{a,b}^n$ are unimodular only if $a=b$. Hence we cannot expect an easy formula for $a\neq b$, also we have already seen that the number of lattice points in $T_{a,a}^n$ is given by
\begin{align*}
d\left(\left\lfloor\frac{n}{a}\right\rfloor\right).
\end{align*}\vspace{1mm}
Instead consider $T_{a,b}^n=\{x,y\in\mathbb{R}_+^2: ax+by\leq n\}$ as lying in the hyperplane $z\equiv n$ in $\mathbb{R}^3$. Then
\begin{align*}
\bigcup_{n\geq 0} T_{a,b}^n\cap\mathbb{Z}^3=\mathrm{co}\left(\left(
                                                             \begin{array}{c}
                                                               1 \\
                                                               0 \\
                                                               a \\
                                                             \end{array}
                                                           \right),\left(
                                                             \begin{array}{c}
                                                               0 \\
                                                               1 \\
                                                               b \\
                                                             \end{array}
                                                           \right),\left(
                                                             \begin{array}{c}
                                                               0 \\
                                                               0 \\
                                                               1 \\
                                                             \end{array}
                                                           \right)
\right)\cap\mathbb{Z}^3.
\end{align*}
The latter cone is unimodular and has generating function
\begin{align*}
f(x,y,z)=\frac{1}{(1-xz^a)(1-yz^b)(1-z)}.
\end{align*}
In particular, the number of lattice points in $T_{a,b}^n$ corresponds to the coefficient of $z^n$ of the expansion of $f$ restricted to $x=y=1$. This explains formula (\ref{eq.PowSerL}).\vspace{6mm}

\writesection{Scale Invariance} The condition in Theorem 1 is scale invariant, meaning that for
each real $\lambda>0$ one has
\begin{align*}
\mathcal{N}(a,b)\preceq\mathcal{N}(c,d) \Longleftrightarrow \mathcal{N}(\lambda a,\lambda
b)\preceq\mathcal{N}(\lambda c,\lambda d).
\end{align*}
Unfortunately, this scale invariance does not descend to the generating functions. Thus $g_{a,b}\succeq g_{c,d}$ does
not imply $g_{\lambda a,\lambda b}\succeq g_{\lambda c,\lambda d}$ and it does not make sense to extend
our notion of generating functions to real parameters $a,b$. For rational $a,b,c,d\in\mathbb{Q}_+$ the best one could
do is to choose $N\in\mathbb{N}$ such that $Na,Nb,Nc,Nd$ are integers and then compare
the generating functions $g_{Na,Nb}$ and $g_{Nc,Nd}$.\vspace{1mm}

The embedding condition $g_{a,b}\succeq g_{c,d}$ requires
\begin{align}\label{eq.scaleinvcond}
g_{a,b}(z)\geq g_{c,d}(z)
\end{align}
for all $z\in [0,1)$. But (\ref{eq.scaleinvcond}) is scale invariant, since it is equivalent to
\begin{align*}
\frac{(1-z^c)(1-z^d)}{(1-z^a)(1-z^b)}\geq 1
\end{align*}
and one may substitute $z=w^\lambda$ with $w\in [0,1)$ on the left hand side. Therefore it corresponds to an embedding
obstruction which extends to real parameters $a,b$. The following lemma shows that at least in the case
of embeddings into a ball this obstruction is the volume constraint.

\writelemma\textsl{Let $a,b,c,d\in\mathbb{R}$ be positive, such that $b\leq\min(c,d)$. Then the inequality
\begin{align*}
g_{a,b}(z)\geq g_{c,d}(z)
\end{align*}
holds for all $z\in[0,1)$ if and only if $a$ is chosen such that $ab\leq cd$.}\vspace{1mm}

{\sc Proof:} By scale invariance it suffices to show that under the assumption $b\leq\min(1,c)$ the inequality
\begin{align}\label{eq.ineqzabc}
(1-z)(1-z^c)\geq (1-z^a)(1-z^b)
\end{align}
holds for all $z\in (0,1)$ if and only if $a\leq \frac{c}{b}$.\vspace{1mm}

We first consider the case $c=ab$, such that $b\leq 1\leq a$. Then we have
\begin{align*}
ab\leq\min(a,ab+1)\leq\max(a,ab+1)\leq a+b.
\end{align*}
The function $f(x)=z^x$ is convex and monotone decreasing for fixed $z\in(0,1)$ and $x\in(0,\infty)$.
Hence the segment from $(ab,z^{ab})$ to $(a+b,z^{a+b})$ lies above the segment from $(a,z^a)$ to $(ab+1,z^{ab+1})$.
Comparing the heights of intersection of these segments with the horizontal line $x=\frac{b(ab+1)+a}{b+1}$ yields the
estimate
\begin{align*}
\frac{b}{b+1}z^{ab+1}+\frac{1}{b+1}z^a\leq\frac{b}{b+1}z^{ab}+\frac{1}{b+1}z^{a+b}.
\end{align*}
Considering the function $F:[1,\infty)\rightarrow\mathbb{R}$,
\begin{align*}
F(a)=z^{ab+1}+z^a+z^b-z^{ab}-z^{a+b}-z
\end{align*}
for fixed $z\in(0,1)$ and $b\leq 1$, the previous inequality implies that $f$ is monotone increasing in
$a$. Consequently, $F(a)\geq F(1)=0$. This tells us that (\ref{eq.ineqzabc}) holds for all $z\in (0,1)$ if $c=ab$. Since
increasing $c$ only increases the left hand side of (\ref{eq.ineqzabc}), we have shown that this inequality is satisfied
for all $z\in (0,1)$ if $c\geq ab$. \vspace{1mm}

Now we fix any $0<\lambda<1$ and consider the case $c=\lambda ab$. Let
\begin{align*}
C:=\frac{\lambda^2ab^2+a+b}{\lambda b+1}>\frac{1+b}{\lambda b+1}>1.
\end{align*}
Choose $\delta>0$ small enough, such that
\begin{align*}
z^{C-\lambda b}\geq -\frac{(a+b)^2}{4(1-\lambda)b}\log z
\end{align*}
holds for $z\in(1-\delta,1)$. Using this and the convexity and monotonicity of the function $f$, we obtain for $1\leq
\tau\leq a$
\begin{align*}
\frac{\lambda b}{\lambda b+1}z^{\lambda\tau b+1}+&\frac{1}{\lambda b+1}z^{\tau}\geq
f\left(\frac{\lambda^2\tau b^2+\tau+\lambda b}{\lambda b+1}\right)\geq f\left(\frac{\lambda^2\tau b^2+\tau+b}{\lambda
b+1}\right)-\frac{(1-\lambda)b}{\lambda b+1}f'\left(\frac{\lambda^2\tau b^2+\tau+b}{\lambda
b+1}\right) \\ &> f\left(\frac{\lambda^2\tau b^2+\tau+b}{\lambda
b+1}\right)-\frac{(1-\lambda)b}{2}f'(C)=f\left(\frac{\lambda^2\tau b^2+\tau+b}{\lambda
b+1}\right)-\frac{(1-\lambda)b}{2}z^C\log z\\ &\geq f\left(\frac{\lambda^2\tau b^2+\tau+b}{\lambda
b+1}\right)+\frac{(a+b)^2}{8}z^{\lambda b}(\log z)^2=f\left(\frac{\lambda^2\tau b^2+\tau+b}{\lambda
b+1}\right)+\frac{(a+b)^2}{8}f''(\lambda b).
\end{align*}
We now apply the inequality
\begin{align*}
\left|\mu f(x)+(1-\mu)f(y)-f\left(\mu x+(1-\mu)y\right)\right|\leq\frac{|x-y|^2}{8}\cdot\max_{\xi\in[x,y]} f''(\xi)
\end{align*}
with $\mu=\frac{\lambda b}{\lambda b+1}$ to conclude
\begin{align*}
\frac{\lambda b}{\lambda b+1}z^{\lambda\tau b+1}+\frac{1}{\lambda b+1}z^{\tau}> & f\left(\mu(\lambda\tau
b)+(1-\mu)(\tau+b)\right)+\frac{|\lambda\tau b-(\tau+b)|^2}{8}\cdot\max_{\xi\in[\lambda\tau
b,\tau+b]} f''(\xi)\\ \geq & \mu f\left(\lambda\tau b\right)+(1-\mu)f\left(\tau+b\right)=\frac{\lambda b}{\lambda
b+1}z^{\lambda\tau b}+\frac{1}{\lambda b+1}z^{\tau+b}
\end{align*}
for $1\leq\tau\leq a$ and $z\in(1-\delta,1)$. Consequently, the function $F_\lambda:[1,a]\rightarrow\mathbb{R}$ defined
by
\begin{align*}
F_\lambda(\tau)=z^{\lambda\tau b+1}+z^{\tau}+z^b-z^{\lambda\tau b}-z^{\tau+b}-z
\end{align*}
is monotone decreasing for $z\in(1-\delta,1)$. Hence for these values of $z$ we have
\begin{align*}
F_\lambda(\tau)\leq F_\lambda(1)=(1-z)(z^b-z^{\lambda b})<0.
\end{align*}
This shows that (\ref{eq.ineqzabc}) is violated for $c=\lambda ab$ with $0<\lambda<1$. \hfill
$\Box$\vspace{3mm}

\paragraph{Acknowledgements.} I warmly thank Dusa McDuff for introducing me into the topic at Edifest 2010 and giving
helpful comments. I also thank Felix Schlenk and Matthias Schwarz for making it possible for me to participate at this
conference. Finally, I thank the Max Planck Institute for Mathematics in the Sciences for support and providing a
pleasant environment to do this research.

\end{document}